\documentclass{article}

\usepackage[utf8]{inputenc}
\usepackage[active]{srcltx}      

\usepackage{amssymb}
\usepackage{amsmath, amsthm, amsfonts}
\usepackage{mathrsfs}
\usepackage{frcursive}
\usepackage{color}

\newcommand\R{{\mathbb R}}

\newcommand\C{{\mathbb C}}

\def\AA{{\cal A}}
\def\BB{{\cal B}}

\def\LL{{\cal L}}

\def\XX{X}

\newtheorem{theo}{Theorem}

\newtheorem{hyp}[theo]{Hypothesis}
\newtheorem{lem}[theo]{Lemma}

\newtheorem{rem}[theo]{Remark}

\newcommand{\beqn}{\begin{equation}}
\newcommand{\eeqn}{\end{equation}}
\newcommand{\bear}{\begin{eqnarray}}
\newcommand{\eear}{\end{eqnarray}}
\newcommand{\bean}{\begin{eqnarray*}}
\newcommand{\eean}{\end{eqnarray*}}

\DeclareMathOperator{\sign}{sign}

\author{
  \textsc{María J. Cáceres}\footnote{Departamento de Matemática
    Aplicada, Universidad de Granada, E18071 Granada, Spain. Email:
    \texttt{caceresg@ugr.es}}
  \and
  \textsc{José A. Cañizo}\footnote{Departament de Matemàtiques,
    Universitat Autònoma de Barcelona, 08193 Bellaterra, Spain. Email:
    \texttt{canizo@mat.uab.es}}
  \and
  \textsc{Stéphane Mischler}\footnote{IUF and CEREMADE,
    Univ. Paris-Dauphine, Place du Maréchal de Lattre de Tassigny,
    75775 Paris CEDEX 16, France. Email:
    \texttt{mischler@ceremade.dauphine.fr}}
}

\date{December 20, 2010}

\begin{document}

\title{Rate of convergence to self-similarity for the fragmentation
  equation in $L^1$ spaces}

\maketitle

\begin{abstract}
  In a recent result by the authors \cite{CCM11} it was proved that
  solutions of the self-similar fragmentation equation converge to
  equilibrium exponentially fast. This was done by showing a spectral
  gap in weighted $L^2$ spaces of the operator defining the time
  evolution. In the present work we prove that there is also a
  spectral gap in weighted $L^1$ spaces, thus extending exponential
  convergence to a larger set of initial conditions. The main tool is
  an extension result in \cite{GMM**}.
\end{abstract}

\section{Introduction}

In a recent paper \cite{CCM11} we have studied the speed of
convergence to equilibrium for solutions of equations involving the
fragmentation operator and first-order differential terms. In this
paper we will focus on the case of self-similar fragmentation given by
\begin{subequations}
  \label{eq:ssf}
  \begin{gather}
    \label{eq:ssf-eq}
    \partial_t g_t(x) 
    = - x \partial_x g_t(x) - 2 g_t(x) + \LL g_t (x)
    \\
    \label{eq:drift-frag-initial}
    g_0(x) = g_{in}(x) \qquad (x > 0).
  \end{gather}
\end{subequations}
Here the unknown is a function $g_t(x)$ depending on time $t \geq 0$
and on size $x > 0$, which represents a density of units (usually
particles, cells or polymers) of size $x$ at time $t$, and $g_{in}$ is
an initial condition. The \emph{fragmentation operator} $\LL$ acts on
a function $g = g(x)$ as
\begin{equation}
  \label{eq:frag-op}
  \LL g(x) :=  \LL_+g(x)  - B(x) g(x),
\end{equation}
where the positive part $\mathcal{L}_+$ is given by
\begin{equation}
  \label{eq:L+}
  \LL_+ g(x) :=  \int_x^\infty b(y,x) g(y) \,dy.
\end{equation}
The coefficient $b(y,x)$, defined for $y > x > 0$, is the
\emph{fragmentation coefficient}, and $B(x)$ is the \emph{total
  fragmentation rate} of particles of size $x > 0$. It is obtained
from $b$ through
\begin{equation}
  \label{eq:B}
  B(x) := \int_0^x \frac{y}{x} \, b(x,y) \,dy
  \qquad (x > 0).
\end{equation}
We refer to
\cite{CCM11,MD86,MR2270822,DG10,citeulike:4788671,citeulike:4308171}
for a motivation of \eqref{eq:ssf} in several applications and a
general survey of the mathematical literature related to it.

We call $T$ the operator on the right hand side of
\eqref{eq:ssf-eq}, this is,
\begin{equation}
  \label{eq:def-T}
  T g (x) := - x \partial_x g(x) - 2 g(x)
  + \LL g (x)
  \qquad (x > 0),
\end{equation}
acting on a (sufficiently regular) function $g$ defined on
$(0,+\infty)$. Notice that, even though $g$ is a one-variable
function, we still denote its derivative as $\partial_x g$ in order to
be consistent with the notation in 
 \eqref{eq:ssf}. The results in
\cite{CCM11} show that $T$ has a spectral gap in the space $L^2(x\,G^{-1})$,
 where $G$ is the unique stationary solution of
 \eqref{eq:ssf} with $\int x \, G = 1$. In the rest of this paper $G$
 will represent this solution, called the \emph{self-similar
   profile}. Proofs of existence of the profile $G$ and some estimates
 are given in \cite{MR2114413,citeulike:4788671,DG10}, and additional
 bounds are given in \cite{CCM11}.

 The main result in \cite{CCM11} is a study of the long time behavior
 of \eqref{eq:ssf}: by means of an inequality relating the quadratic
 entropy and its dissipation rate, exponential convergence is obtained
 in $L^2(x \, G^{-1})$. Using the results in \cite{GMM**} this is
 further extended to the space $L^2(x + x^k)$ for a sufficiently large
 exponent $k$. In this way one obtains convergence in a strong norm,
 but correspondingly has to impose more on the initial condition than
 just having finite mass.

 The purpose of this work is to prove that $T$ has a spectral gap in
 the larger spaces $L^1(x^m + x^M)$, where $1/2 < m < 1 < M$ are
 suitable exponents. This extension is an example of application of
 the results in \cite{GMM**}. The interest of this concerning the
 asymptotic behavior of \eqref{eq:ssf} is that it shows exponential
 convergence is valid for more general initial conditions (any
 function in $L^1(x^m + x^M)$).

\paragraph{Assumptions on the fragmentation coefficient}

In order to use the results in \cite{CCM11} we will make the following
hypotheses on the fragmentation coefficient $b$:
\begin{hyp}
  \label{hyp:b-basic}
  For all $x > 0$, $b(x,\cdot)$ is a nonnegative measure on the
  interval $[0,x]$. Also, for all $\psi \in
  \mathcal{C}_0([0,+\infty))$, the function $x \mapsto \int_{[0,x]}
  b(x,y) \psi(y) \,dy$ is measurable.
\end{hyp}

\begin{hyp}
  \label{hyp:number-of-pieces}
  There exists $\kappa > 1$ such that
  \begin{equation}
    \label{eq:number-of-pieces}
    \int_0^x b(x,y)\,dy = \kappa B(x)
    \qquad (x > 0).
  \end{equation}
\end{hyp}

\begin{hyp}
  \label{hyp:p-bounded}
  There exists $0 <B_m < B_M$ satisfying
  \begin{equation}
    \label{eq:hypotheses-simple}
    2 B_m\, x^{\gamma-1} \leq b(x,y) \leq 2 B_M\, x^{\gamma-1}
    \qquad (0 < y < x)
  \end{equation}
  for some $0 < \gamma < 2$.
\end{hyp}
This implies the following useful bound, as remarked in
\cite[Corollary 6.4]{CCM11}:
\begin{lem}
  \label{lem:b-moment-bounds}
  Consider a fragmentation coefficient $b$ satisfying Hypotheses
  \ref{hyp:b-basic}--\ref{hyp:p-bounded}. There exists a strictly
  decreasing function $k \mapsto p_k$ for $k \geq 0$ with $\lim_{k \to
    +\infty} p_k = 0$,
  \begin{equation}
    \label{eq:Pk-range}
    p_k > 1 \text{ for } k \in [0,1),
    \quad p_1 = 1,
    \quad 0 < p_k < 1 \text{ for } k > 1,
  \end{equation}
  and such that
  \begin{equation}
    \label{eq:b-moments}
    \int_0^x y^k b(x,y) \,dy
    \leq
    p_k\, x^k B(x)
    \qquad (x > 0, \ k > 0).
  \end{equation}
\end{lem}

\paragraph{Main results}

The main result of the present work is a spectral gap of $T$
 on weighted $L^1$ spaces.

\begin{theo}
  \label{thm:main}
  Assume hypotheses \ref{hyp:b-basic}--\ref{hyp:p-bounded}. For any
  $1/2 < m < 1$ there exists $1 < M < 2$ such that the operator
  (\ref{eq:def-T}) has a spectral gap in the space $\XX := L^1(x^m +
  x^M)$. More precisely, there exists $\alpha > 0$ and a constant $C
  \geq 1$ such that, for all $g_{in} \in \XX$ with $\int x \, g_{in} =
  1$
  \begin{equation*}
    \|g_t - G\|_{\XX}
    \leq C \,e^{-\alpha t}\, \|g_{in} - G\|_{\XX}
    \qquad (t \geq 0).
  \end{equation*}
\end{theo}

\section{Preliminaries}
\label{sec:prelims}

In this section we gather some known results from previous works.

\subsection{Previous results on the spectral gap of $T$}

A result like Theorem \ref{thm:main} was proved in \cite{CCM11}, but
in the $L^2$ space with weight $x\, G^{-1}$. This is summarized in the
following theorem:

\begin{theo}[\cite{CCM11}]
  \label{thm:L2-spectral-gap}
  Assume Hypotheses \ref{hyp:b-basic}--\ref{hyp:p-bounded}, and
  consider $G$ the self-similar profile with $\int x\,G = 1$. The
  operator $T$ given by \eqref{eq:def-T} has a spectral gap in the
  space $H = L^2(x \,G^{-1})$.

  More precisely, there exists $\beta > 0$ such that for any $g_{in}
  \in H$ with $\int x\,g = 1$ the solution $g \in C([0,\infty);
  L^1(x\,dx))$ to equation \eqref{eq:ssf} satisfies
  \begin{equation*}
    \qquad \| g_t - G \|_H \le e^{- \beta \, t}
    \, \| g_{in} - G \|_H
    \qquad (t \geq 0).
  \end{equation*}
\end{theo}

\subsection{Bounds for the self-similar profile}

We recall the following result from \cite[Theorem 3.1]{CCM11}:
\begin{theo}
  \label{thm:ssf-bounds}
  Assume Hypotheses \ref{hyp:b-basic}--\ref{hyp:p-bounded} on the
  fragmentation coefficient $b$, and call $\Lambda(x) := \int_0^x
  \frac{B(s)}{s} \,ds$. Let $G$ be the self-similar profile with $\int
  x\,G = 1$.

  For any $\delta > 0$ and any $a \in (0, B_m/B_M)$, $a' \in
  (1,+\infty)$ there exist constants $C' = C'(a',\delta)$, $C = C(a) >
  0$ such that
  \begin{equation} 
    \label{eq:ssf-bounds}
    C'\, e^{-a' \Lambda(x)} \leq
    G(x)
    \leq  C\, e^{-a\, \Lambda(x)}
    \quad \text{ for } x > 0.
  \end{equation}
\end{theo}

\begin{rem}
  In the case $b(x,y)=2\, x^{\gamma-1}$ (so $B(x)=x^\gamma$), the
  profile $G$ has the explicit expression $G(x) = e^
  {-\frac{x^\gamma}{\gamma}} $ for $\gamma>0$. This motivates the
  choice of $e^{-a\, \Lambda(x)}$ as functions for comparison. For a
  general $b(x,y)$ no explicit form is available.
\end{rem}

\begin{proof}
  Everything but the lower bound of $G$ for small $x$ is proved in
  \cite[Section 3]{CCM11}. For the lower bound, we calculate as
  follows:
  \begin{equation}
    \label{eq:dxK}
    \partial_x \left(x^2\, e^{\Lambda(x)} G(x) \right)
    = x\, e^{\Lambda(x)} \int_x^\infty b(y,x) \,G(y) \, dy
    \qquad (x > 0),
  \end{equation}
  which implies that $x^2\, e^{\Lambda(x)} G(x)$ is a nondecreasing
  function. Hence, it must have a limit as $x \to 0$, and this limit
  must be $0$ since we know $x\,G(x)$ is integrable. Then, integrating
  \eqref{eq:dxK}, and for $0 < z < 1$,
  \begin{multline}
    \label{eq:dxK-2}
    z^2\, e^{\Lambda(z)} G(z)
    = \int_0^z x\, e^{\Lambda(x)} \int_x^\infty b(y,x) \,G(y) \, dy \,dx
    \\
    = \int_0^\infty \,G(y) \int_0^{\min\{z,y\}}
    b(y,x) x\, e^{\Lambda(x)} \, dx \,dy
    \\
    \geq 2 B_m \int_0^\infty y^{\gamma-1}\,G(y) \int_0^{\min\{z,y\}}
    x \, dx \,dy
    \\
    =
    B_m \int_0^\infty y^{\gamma-1}\,G(y) (\min\{z,y\})^2 \,dy
    \\
    \geq
    B_m z^2 \int_z^\infty y^{\gamma-1}\,G(y) \,dy
    \\
    \geq
    B_m z^2 \int_1^\infty y^{\gamma-1}\,G(y) \,dy
    = C z^2
    \qquad (0 < z < 1).
  \end{multline}
  Notice that the number $\int_1^\infty y^{\gamma-1}\,G(y) \,dy$ is
  strictly positive, as the profile $G$ is strictly positive
  everywhere (see \cite{DG10,MR2114413,CCM11}). This proves the lower
  bound on $G(x)$ for $0 < x < 1$, and completes the proof.
\end{proof}

\subsection{A general spectral gap extension result}

Our proof is based on the following result from \cite{GMM**}, which
was already used in \cite{CCM11} for an extension to an $L^2$ space
with a polynomial weight:
\begin{theo}
  \label{thm:GMM**}
  Consider a Hilbert space $H$ and a Banach space $\XX$ (both over the
  field $\C$ of complex numbers) such that $H \subset \XX$ and $H$ is
  dense in $\XX$. Consider two unbounded closed operators with dense
  domain $T$ on $H$, $\Lambda$ on $\XX$ such that $\Lambda_{|H}= T$.
  On $H$ assume that
  \begin{enumerate}
  \item \label{it:nullspace} There is $G \in H$ such that $T \, G = 0$
    with $\| G \|_H = 1$;
    
  \item \label{it:invariant-quantity} Defining $\psi(f) := \langle f,G
    \rangle_H$, the space $ H_0 := \{ f \in H; \,\, \psi(f) = 0 \}
    $ is invariant under the action of $T$.

  \item \label{it:Ldissipative} $T- a$ is dissipative on $H_0$ for
    some $a < 0$, in the sense that
    $$
    \forall \, g \in D(T) \cap H_0 \qquad ((T-a) \, g, g)_H \le 0,
    $$
    where $D(T)$ denotes the domain of $T$ in $H$.

  \item \label{it:Lsemigroup} $T$ generates a semigroup $e^{t \, T}$ on $H$;

    \smallskip\noindent Assume furthermore on $\XX$ that
    
  \item \label{it:A+B}
    there exists a continuous linear form $\Psi : \XX \to \R$ such
    that $\Psi_{|H} = \psi$;
    
    \smallskip\noindent and $\Lambda$ decomposes as $\Lambda = \AA + \BB$
    with

  \item \label{it:Abounded} $\AA$ is a bounded operator from $\XX$ to $H$;

  \item \label{it:Bdissipative} $\BB$ is a closed unbounded operator
    on $\XX$ (with same domain as $D(\Lambda)$ the domain of
    $\Lambda$) and the semigroup $e^{t\BB}$ it generates satisfies,
    for some constant $C \geq 1$, that
    \begin{equation}
      \label{eq:B-dissipative}
      \forall t \geq 0, \
      \forall g \in \XX  \text{ with } \Psi(g)=0,
      \quad
      \|e^{t\BB} g\|_{\XX} \leq C \|g\|_{\XX}\, e^{a t},
    \end{equation}
    where $a < 0$ is the one from point \ref{it:Ldissipative}.
  \end{enumerate}
  \smallskip\noindent
  Then, for any $a' \in (a,0)$ there exists $C_{a'} \ge 1$ such that
  $$
  \forall \, t \ge 0,\
  \forall g \in \XX,
  \qquad
  \| e^{t \Lambda} \, g - \Psi(g) \, G \|_\XX
  \le C_{a'} \, \| g - \Psi(g) \, G \|_\XX \, e^{a' t}.
  $$
\end{theo}

\section{Proof of the main theorem}


The proof consists is an application of Theorem \ref{thm:GMM**}. For
this, we consider the Hilbert space $H := L^2(x\,G^{-1}(x))$, where
$G$ is the unique self-similar profile with $\int G = 1$, and define
$\psi(g) := \int x g$. Due to our previous results \cite{CCM11} we
know that $T$ and $\psi$ satisfy points
\ref{it:nullspace}--\ref{it:Lsemigroup} of Theorem \ref{thm:GMM**}.

As the larger space we take $\XX = L^1(x^m + x^M)$, with $1/2 < m < 1
< M$, to be precised later. Observe that, due to the bounds on $G$
from Theorem \ref{thm:ssf-bounds},
\begin{multline*}
  \| g \|_{\XX}
  =
  \int_0^\infty (x^m + x^M) |g(x)| \,dx
  \\
  \leq
  \left(\int_0^\infty g(x)^2 \frac{x}{G(x)} \,dx \right)^{1/2}
  \left( \int_0^\infty (x^{m-\frac{1}{2}} + x^{M-\frac{1}{2}})^2
    G(x) \,dx \right)^{1/2}
  = C \|g\|_{H},
\end{multline*}
and hence $H \subseteq \XX$. Similarly,
\begin{equation*}
  \int_0^\infty x |g(x)| \,dx \leq
  \int_0^\infty (x^m + x^M) |g(x)| \,dx,
\end{equation*}
which allows us to define $\Psi: \XX \to \R$, $\Psi(g) := \int x g$,
and proves that $\Psi$ is continuous on $\XX$. Obviously $\Psi_{|H} =
\psi$, so point \ref{it:A+B} of Theorem \ref{thm:GMM**} is also
satisfied.

Consider $\Lambda$ the unbounded operator on $\XX$ given by the same
expression \eqref{eq:def-T} (with domain a suitable dense subspace of
$\XX$ which makes $\Lambda$ a closed operator). To prove the remaining
points \ref{it:Abounded} and \ref{it:Bdissipative} we use the
following splitting of $\Lambda$, taking real numbers $0 < \delta < R$
to be chosen later:
\begin{gather}
  \label{eq:def-A-ssf}
  \begin{split}
    \mathcal{A}g (x) := \LL^{+,s} g(x)
    &:= \int_x^\infty b_{R,\delta}
    (y,x) \, g(y) \, dy
    \\
    &= {\bf 1}_{x \le R} \int_x^\infty {\bf
      1}_{y \ge \delta} \, b(y,x) g(y) \, dy,
  \end{split}
  \\
  \label{eq:Lambda=A+B}
  \Lambda = \mathcal{A} + \mathcal{B},
  \qquad \mathcal{B} g := \Lambda g - \mathcal{A} g,
\end{gather}
where we denote $b_{R,\delta} (x,y) := b(x,y) \, {\bf 1}_{x \ge
  \delta} \, {\bf 1}_{y \le R}$. We define
\bean \LL^{+,r} g &:=&
\LL^+ g - \LL^{+,s} g
\\
&=& \int_x^\infty b(y,x) \, ( 1 - {\bf 1}_{y \ge \delta} \, {\bf 1}_{x
  \le R}) \, g(y) \, dy
\\
&=& \int_x^\infty b(y,x) \, {\bf 1}_{y \le \delta} \, g(y) \, dy +
\int_x^\infty b(y,x) \, {\bf 1}_{y \ge \delta} \, {\bf 1}_{x \ge R} \,
g(y) \, dy
\\
&=:& \LL_1^{+,r} g + \LL_2^{+,r} g
\eean
so we may write $\mathcal{B}$ as
\begin{equation}
  \label{eq:def-B-2}
  \mathcal{B}g
  = -2 g - x \partial_x g - B g + \LL_1^{+,r} g + \LL_2^{+,r} g.
\end{equation}

First, let us prove that $\mathcal{A}$ is bounded from $\XX$ to
$H$. We compute
\begin{align*}
  \|\mathcal{A} g \|_{H}^2
  &=
  \int_0^\infty x \, (\LL^{+,s} g )^2 \, G(x)^{-1} \, dx
  \\
  &\leq (2 B_M)^2  \left(
    \sup_{[0,R]} x \, G(x)^{-1}\right) \int_0^R \, \left(
    \int_{\max(x,\delta)}^\infty y^{\gamma-1} g(y) \, dy \right)^2 \,
  dx
  \\
  &\le C_R \left( \int_\delta^\infty y^{\gamma-1} \, g(y) \, dy
  \right)^2
  \leq
  C_{R,\delta} \left( \int_0^\infty y \, g(y) \, dy \right)^2
  \leq C_{R,\delta} \, \|g\|_{\XX}^2,
\end{align*}
which shows $\mathcal{A}:\XX \to H$ is a bounded operator. Notice that
we have used here the lower bound $G(x) \geq C x$ for $x$ small,
proved in Theorem \ref{thm:ssf-bounds}.

Then, let us prove that one can choose $0 < \delta < R$ appropriately
so that $\mathcal{B}$ satisfies point \ref{it:Bdissipative} of Theorem
\ref{thm:GMM**} for some $a < 0$. It is enough to prove that, for $g$
a real function in the domain of $\Lambda$ (the same as the domain of
$\mathcal{B}$),
\begin{equation}
  \label{eq:B-a-dissipative-aim}
  \int_0^\infty  \sign(g(x)) \,\mathcal{B} g(x)\, (x^m + x^M) \,dx
  \leq
  a \|g\|_{\XX},
\end{equation}
since then one can obtain \eqref{eq:B-dissipative} with $C=1$ by
considering the time derivative of the $L^1$ norm of $e^{t\BB}g$. If
we have this for any real $g$, it is easy to show it also holds for a
complex $g$ and some constant $C \geq 1$.
So, we take $g$ real and in the domain of $\Lambda$, and calculate as
follows for any $k > 0$, using \eqref{eq:def-B-2}:
\begin{multline}
  \label{eq:B-bound-1}
  \int_0^\infty  \sign(g(x)) \,\mathcal{B} g(x)\, x^k \,dx
  \leq
  (k-1) \int_0^\infty x^k \, |g| \, dx
  \\
  - \int_0^\infty B(x) x^{k} \, |g| \, dx
  + \int_0^\infty |\LL_1^{+,r} g| x^k \,dx
  + \int_0^\infty |\LL_2^{+,r} g| x^k \,dx,
\end{multline}
where the first term is obtained from the terms $-2g-\partial_x g$
through an integration by parts. We give separately some bounds on the
last two terms in \eqref{eq:B-bound-1}. On one hand, we have
\begin{equation}
  \label{eq:L_1-bound}
\begin{split}
  \int_0^\infty |\LL_1^{+,r} g| x^k \,dx
  & \le \int_0^\infty x^k
  \int_x^\infty b(y,x) \, {\bf 1}_{y \le \delta} \, |g(y)| \, dy \, dx
  \\
  & \le \int_0^\delta |g(y)| \, \Bigl( \int_0^y x^k \, b(y,x) \, dx
  \Bigr) \, dy
  \\
  & \le 2 B_M \int_0^\delta |g(y)| B(y) y^k  \, dy
  \\
  & \le p_k B_m \delta^\gamma \int_0^\delta |g(y)| y^k  \, dy,
  \end{split}
\end{equation}
where we have used \eqref{eq:b-moments}. On the other hand, and again
due to \eqref{eq:b-moments},
\begin{equation}
  \label{eq:L_2-bound}
  \begin{split}
    \int_0^\infty |\LL_2^{+,r} g| x^k \,dx & \le \int_0^\infty x^k
    \int_x^\infty b(y,x) \, {\bf 1}_{x \ge R} \, {\bf 1}_{y \ge
      \delta} \, |g(y)| \, dy \, dx
    \\
    & \le \int_0^\infty x^k \int_x^\infty b(y,x) \, {\bf 1}_{x \ge R}
    \, {\bf 1}_{y \ge R} \, |g(y)| \, dy \, dx
    \\
    & \le \int_R^\infty |g(y)| \, \Bigl( \int_R^y x^k \, b(y,x) \, dx
    \Bigr) \, dy
    \\
    & \le p_k \int_R^\infty |g(y)| y^k B(y) \, dy.
  \end{split}
\end{equation}
Hence, from \eqref{eq:B-bound-1} and the bounds
\eqref{eq:L_1-bound}--\eqref{eq:L_2-bound} we obtain
\begin{multline}
  \label{eq:B-bound-2}
  \int_0^\infty \mathcal{B}g(x) \sign(g(x)) (x^m + x^M) \,dx
  \\
  \leq
  (m-1) \int_0^\infty x^m \, |g| \, dx
  +
  (M-1) \int_0^\infty x^M \, |g| \, dx
  \\
  - \int_0^\infty B(x) (x^m + x^M) \, |g| \, dx
  \\
  + p_m B_m \delta^\gamma \int_0^\delta x^{m} |g(x)| \, dx
  + p_m \int_R^\infty x^{m} B(x) \, |g(x)| \, dx
  \\
  + p_M B_m \delta^\gamma \int_0^\delta x^{M} |g(x)| \, dx
  + p_M  \int_R^\infty x^{M} B(x) \, |g(x)| \, dx.
\end{multline}
We have to choose $1/2 < m < 1 < M < 2$ so that this is bounded by $-
C \|g\|_{\XX}$ for some positive constant $C$. First, fix any $m$ with
$1/2 < m < 1$, and take $0 < \delta < 1$ small enough such that
\begin{equation*}
  p_m B_m \delta^\gamma < \frac{1-m}{4},
  \quad B_m \delta^\gamma < \frac{1-m}{4}.
\end{equation*}
(Which can be done due to $\gamma > 0$.) Then, as $p_M < 1$ and $x^M <
x^m$ for $x < \delta < 1$,
\begin{multline}
  \label{eq:B-bound-2.5}
  \int_0^\infty \mathcal{B}g(x) \sign(g(x)) (x^m + x^M) \,dx
  \\
  \leq
  -\frac{1-m}{2} \int_0^\infty x^m \, |g| \, dx
  +
  (M-1) \int_0^\infty x^M \, |g| \, dx
  \\
  - \int_0^\infty B(x) (x^m + x^M) \, |g| \, dx
  \\
  + p_m \int_R^\infty x^{m} B(x) \, |g(x)| \, dx
  + p_M  \int_R^\infty x^{M} B(x) \, |g(x)| \, dx.
\end{multline}
Now, take $R_0 > 0$ such that $B(x) > 2 > M$ for $x \geq R_0$. Then,
choose $1 < M < 2$ such that $(M-1) x^M < \frac{1-m}{4} x^m$ for $0 <
x < R_0$. Then whatever $R$ is we have from \eqref{eq:B-bound-2}:
\begin{multline}
  \label{eq:B-bound-3}
  \int_0^\infty \mathcal{B}g(x) \sign(g(x)) (x^m + x^M) \,dx
  \\
  \leq
  -\frac{1-m}{4} \int_0^{R_0} x^m \, |g| \, dx
  - \int_{R_0}^R x^M \, |g| \, dx
  \\
  - \int_{R}^\infty (B(x) - M + 1) x^M \, |g| \, dx
  \\
  + p_m \int_R^\infty x^{m} B(x) \, |g(x)| \, dx
  + p_M  \int_R^\infty x^{M} B(x) \, |g(x)| \, dx.
\end{multline}
Finally, choose $R > 1$ such that
\begin{equation*}
  -(B(x)(1-p_M) - M + 1) x^M + p_m x^m
  \leq - x^M
  \quad \text{ for } x > R.
\end{equation*}
With this, and continuing from \eqref{eq:B-bound-3},
\begin{multline}
  \label{eq:B-bound-4}
  \int_0^\infty \mathcal{B}g(x) \sign(g(x)) (x^m + x^M) \,dx
  \\
  \leq
  -\frac{1-m}{4} \int_0^{R_0} x^m \, |g| \, dx
  - \int_{R_0}^\infty x^M \, |g| \, dx
  \\
  \leq -C \, \|g\|_{\XX},
\end{multline}
for some number $C=C(m,M,R_0) > 0$. This shows point that
$\mathcal{B}$ is dissipative with constant $-C$, and hence point
\ref{it:Bdissipative} of Theorem \ref{thm:GMM**} holds with $a =
-C$. A direct application of Theorem \ref{thm:GMM**} then proves our
result, Theorem \ref{thm:main}, with $\alpha := \min\{\beta, C\}$
(where $\beta$ is the one appearing in Theorem
\ref{thm:L2-spectral-gap}).

\bigskip {\footnotesize \noindent\textit{Acknowledgments.} The first
  two authors acknowledge support from the project
  MTM2008-06349-C03-03 DGI-MCI (Spain) and the Spanish-French project
  FR2009-0019. The second author is also supported by the project
  2009-SGR-345 from AGAUR-Generalitat de Catalunya. The third author
  acknowledges support from the project ANR-MADCOF.}


\end{document}